\documentclass[12pt,fleqn, a4]{article}

\usepackage[french]{babel}
\usepackage{inputenc}
\usepackage[T1]{fontenc}

\usepackage[mathscr]{euscript}
\usepackage{amsmath}

\usepackage{amssymb}
\usepackage{latexsym}
\usepackage{graphics}

\usepackage{graphics}
\usepackage{color}
\newcommand{\colval}{0.3}
\definecolor{colone}{gray}{\colval}

\newcommand{\dcb}{\begin{array}{lll}}
\newcommand{\dce}{\end{array}}
\newcommand{\ebe}{\begin{enumerate}\setlength{\baselineskip}{13pt}\setlength{\parskip}{5pt}}
\newcommand{\dbe}{\end{enumerate}}

\newcommand{\ibegin}{\begin{itemize}\setlength{\baselineskip}{19pt}\setlength{\parskip}{7pt}}
\newcommand{\iend}{\end{itemize}}
\newcommand{\ok}{\rule{4pt}{6pt}}
\newcommand{\desb}{\begin{description}}
\newcommand{\dese}{\end{description}}

\newtheorem{Thm}{Theorem}[section]
\newtheorem {Cor}[Thm]{Corollary}
\newtheorem {definition}{Definition}[section]
\newtheorem {pro}{Proposition}[Thm]
\newtheorem {Lemma}[Thm]{Lemma}
\newtheorem {rem}{Remark}[section]

\newtheorem {assumption}{Assumption}[section]
\newtheorem {condition}[definition]{Condition}

\newcommand {\bd}{\begin{definition}}
\newcommand {\ed}{\end{definition}}
\newcommand {\bpro}{\begin{pro}}
\newcommand {\epro}{\end{pro}}
\newcommand {\bl}{\begin{Lemma}}
\newcommand {\el}{\end{Lemma}}
\newcommand {\bcor}{\begin{Cor}}
\newcommand {\ecor}{\end{Cor}}
\newcommand {\brem }{\begin{rem} \rm }
\newcommand {\erem }{\end{rem}}
\newcommand{\bethe}{\begin{Thm}}
\newcommand{\ethe}{\end{Thm}}
\newcommand {\bassumption}{\begin{assumption}}
\newcommand {\eassumption}{\end{assumption}}
\newcommand {\bcondition}{\begin{condition}}
\newcommand {\econdition}{\end{condition}}

\newcommand{\transp}{{^\top\!}}

\def \ind{1\!\!1}

\begin{document}

\begin{center}
\Large
Construction of multi-default models with full viability\footnote{a working version}
\end{center}

\begin{center}
Shiqi Song

{\footnotesize Laboratoire Analyse et Probabilit\'es\\
Universit\'e d'Evry Val D'Essonne, France\\
shiqi.song@univ-evry.fr}
\end{center}

\

We have the following idea why a financial market is not complete. In fact, initially any market is fair, smooth and complete, just like the Black-Scholes model. Only over the time, successive default events and crisis deteriorate the market condition. The market persists, but becomes unpredictable and incomplete. A natural way to model this evolution is to begin with a fair information flow $\mathbb{F}$, and then to expand $\mathbb{F}$ successively with random times $\tau_{1},\ldots, \tau_{n}$ (multi-default time model), where the random times $\tau_{1},\ldots, \tau_{n}$ are chosen so that the market with expanded information flow remains viable. 

Notice that there does not exist many classes of random times with which successive expansion can be made with viability. One thinks naturally the class of honest times. But it is not a good idea, because honest times in general destroy the market viability (cf. \cite{FJS}). One can take the class of initial times introduced in \cite{JC} to do the successive expansion as in \cite{pham}. Such an expansion may preserve the viability because of Jacod's criterion. But in this paper, we try to do the successive expansion with a new class of random times, namely the random times of $\natural$-models introduced in \cite{JS2, songmodel}.

\

\section{Notations and vocabulary}\label{vocabulary}

We employ the vocabulary of stochastic calculus as defined in \cite{HWY, Jacodlivre} with the following specifications.

\

\textbf{Probability space and random variables}

A stochastic basis $(\Omega, \mathcal{A},\mathbb{P},\mathbb{F})$ is a quadruplet, where $(\Omega, \mathcal{A},\mathbb{P})$ is a probability space and $\mathbb{F}$ is a filtration of sub-$\sigma$-algebras of $\mathcal{A}$, satisfying the usual conditions. 

The relationships involving random elements are always in the almost sure sense. For a random variable $X$ and a $\sigma$-algebra $\mathcal{F}$, the expression $X\in\mathcal{F}$ means that $X$ is $\mathcal{F}$-measurable. The notation $\mathbf{L}^p(\mathbb{P},\mathcal{F})$ denotes the space of $p$-times $\mathbb{P}$-integrable $\mathcal{F}$-measurable random variables.

\

\textbf{Vectors}

An element $v$ in an Euclidean space $\mathbb{R}^d$ ($d\in\mathbb{N}^*$) is considered as a vertical vector. We denote its transposition by $\transp v$. The components of $v$ will be denoted by $v_h, 1\leq h\leq d$. The vectors are always defined as vertical vectors. Let $u=(u_{1},\ldots,u_{d})$ and $v=(v_{1},\ldots,v_{d})$ be two vectors. We denote by $\transp u v$ the inner product $u_{1}v_{1}+\ldots+u_{d}v_{d}$. We denote by $u\vdots v$ the vector composed of $u_{i}v_{i}, 1\leq i\leq d$.

\

\textbf{The processes}

The jump process of a c\` adl\`ag process $X$ is denoted by $\Delta X$, whilst the jump at time $t\geq 0$ is denoted by $\Delta_tX$. By definition, $\Delta_0X=0$ for any c\` adl\`ag process $X$. When we call a process $A$ a process having finite variation, we assume automatically that $A$ is c\`adl\`ag. We denote then by $\mathsf{d}A$ the (signed) random measure that $A$ generates.

We deal with finite family of real processes $X=(X_h)_{1\leq h\leq d}$. It will be considered as $d$-dimensional vertical vector valued process. The value of a component $X_h$ at time $t\geq 0$ will be denoted by $X_{h,t}$. When $X$ is a semimartingale, we denote by $[X,\transp X]$ the $d\times d$-dimensional matrix valued process whose components are $[X_i,X_j]$ for $1\leq i,j\leq d$.

\

\textbf{The projections}

With respect to a filtration $\mathbb{F}$, the notation ${^{\mathbb{F}\cdot p}}\bullet$ denotes the predictable projection, and the notation $\bullet^{\mathbb{F}\cdot p}$ denotes the predictable dual projection. 

\

\textbf{The martingales and the semimartingales}

Fix a probability $\mathbb{P}$ and a filtration $\mathbb{F}$. For any $(\mathbb{P},\mathbb{F})$ special semimartingale $X$, we can decompose $X$ in the form (see \cite[Theorem 7.25]{HWY}) :$$
\dcb
X=X_0+X^m+X^v,\
X^m=X^c+X^{da}+X^{di},
\dce
$$
where $X^m$ is the martingale part of $X$ and $X^v$ is the drift part of $X$, $X^c$ is the continuous martingale part, $X^{da}$ is the part of compensated sum of accessible jumps, $X^{di}$ is the part of compensated sum of totally inaccessible jumps. We recall that this decomposition of $X$ depends on the reference probability and the reference filtration. We recall that every part of the decomposition of $X$, except $X_0$, is assumed null at $t=0$.

\

\textbf{The stochastic integrals}

In this paper we employ the notion of stochastic integral only about the predictable processes. The stochastic integral are defined as 0 at $t=0$. We use a point \texttt{"}$\centerdot$\texttt{"} to indicate the integrator process in a stochastic integral. For example, the stochastic integral of a real predictable process ${H}$ with respect to a real semimartingale $Y$ is denoted by ${H}\centerdot Y$, while the expression $\transp{K}(\centerdot[X,\transp X]){K}$ denotes the process$$
\int_0^t \sum_{i=1}^k\sum_{j=1}^k{K}_{i,s}{K}_{j,s} \mathsf{d}[X_i,X_j]_s,\ t\geq 0,
$$
where ${K}$ is a $k$-dimensional predictable process and $X$ is a $k$-dimensional semimartingale. The expression $\transp{K}(\centerdot[X,\transp X]){K}$ respects the matrix product rule. The value at $t\geq 0$ of a stochastic integral will be denoted, for example, by $\transp{K}(\centerdot[X,\transp X]){K}_t$.

The notion of the stochastic integral with respect to a multi-dimensional local martingale $X$ follows \cite{Jacodlivre}. We say that a (multi-dimensional) $\mathbb{F}$ predictable process is integrable with respect to $X$ under the probability $\mathbb{P}$ in the filtration $\mathbb{F}$, if the non decreasing process $\sqrt{\transp{H}(\centerdot[X,\transp X]){H}}$ is $(\mathbb{P},\mathbb{F})$ locally integrable. For such an integrable process ${H}$, the stochastic integral $\transp{H}\centerdot X$ is well-defined and the bracket process of $\transp{H}\centerdot X$ can be computed using \cite[Remarque(4.36) and Proposition(4.68)]{Jacodlivre}. Note that two different predictable processes may produce the same stochastic integral with respect to $X$. In this case, we say that they are in the same equivalent class (related to $X$).

The notion of multi-dimensional stochastic integral is extended to semimartingales. We refer to \cite{JacShi} for details.

\

\section{Full viability}

\subsection{Definition}

A financial market is modeled by a triplet $(\mathbb{P},\mathbb{F},S)$ of a probability measure $\mathbb{P}$ on a measurable space $(\Omega,\mathcal{A})$, of an information flow $\mathbb{F}=(\mathcal{F}_t)_{t\in\mathbb{R}_+}$ (a filtration of sub-$\sigma$-algebra in $\mathcal{A}$), and of an $\mathbb{F}$ asset process $S$ (a multi-dimensional $\mathbb{F}$ special semimartingale with strictly positive components). The notion of viability has been defined in \cite{HK1979} for a general economy. This notion is then used more specifically to signify that the utility maximization problems have solutions in \cite{choulli3, K2010nu, Kar, K2012, K2011se, loew}. The viability is closely linked to the absences of arbitrage opportunity (of some kind) as explained in \cite{loew, KC2010} so that the word sometimes is employed to signify no-arbitrage condition. In this paper, this notion will be involved in a setting of information flow expansion. Let $\mathbb{G}=(\mathcal{G}_t)_{t\geq 0}$ be a second filtration. We say that $\mathbb{G}$ is an expansion (or an enlargement) of the filtration $\mathbb{F}$, if $\mathcal{F}_t\subset\mathcal{G}_t$, and then, we write $\mathbb{F}\subset \mathbb{G}$. 

\bd\label{scdef}
Let $T>0$ be an $\mathbb{F}$ stopping time. We call a strictly positive $\mathbb{F}$ adapted real process $Y$ with $Y_0=1$, a local martingale deflator on the time horizon $[0,T]$ for a (multi-dimensional) $(\mathbb{P},\mathbb{F})$ special semimartingale $S$, if the processes $Y$ and $Y S$ are $(\mathbb{P},\mathbb{F})$ local martingales on $[0,T]$. The same notion can be defined for the filtration $\mathbb{G}$.
\ed

We recall that the existence of local martingale deflators is equivalent to the no-arbitrage conditions \texttt{NUPBR} and \texttt{NA1} (cf. \cite{takaoka, song-takaoka}). We know that, when the no-arbitrage condition \texttt{NUPBR} is satisfied, the market is viable, and vice versa (cf. \cite{KC2010}). For this reason, we introduce the following definition.

\bd\label{fullviabilitydef}
(\textbf{Full viability} on $[0,T]$ for the expansion $\mathbb{F}\subset \mathbb{G}$) Let $T$ be a $\mathbb{G}$ stopping time. We say that the expansion $\mathbb{F}\subset \mathbb{G}$ is fully viable on $[0,T]$ under $\mathbb{P}$, if, for any $\mathbb{F}$ asset process $S$ possessing a $(\mathbb{P},\mathbb{F})$ deflator, the process $S$ possesses a deflator in the expanded market environment $(\mathbb{P},\mathbb{G})$ on the time horizon $[0,T]$. 
\ed

\brem\label{locallyboundedM}
The full viability implies Hypothesis $(H')$ on $[0,T]$.
\erem

\

\subsection{Enlargements of filtrations and Hypothesis$(H')$ }

Consider the two filtrations $\mathbb{F}=(\mathcal{F}_t)_{t\geq 0}$ and $\mathbb{G}=(\mathcal{G}_t)_{t\geq 0}$ such that $\mathcal{F}_t\subset\mathcal{G}_t$. Let $T$ be a $\mathbb{G}$ stopping time. We introduce the Hypothesis$(H')$ (cf. \cite{CJY, J, JY, Pr, mansuyYor}):   

\bd
(\textbf{Hypothesis$(H')$} on the time horizon $[0,T]$) We say that {Hypothesis$(H')$} holds for the expansion pair $\mathbb{F}\subset \mathbb{G}$ on the time horizon $[0,T]$ under the probability $\mathbb{P}$, if all $(\mathbb{P},\mathbb{F})$ local martingale is a $(\mathbb{P},\mathbb{G})$ semimartingale on $[0,T]$.
\ed

Whenever Hypothesis$(H')$ holds, the associated drift operator can be defined (cf. \cite{song-mrp-drift}).

\bl\label{linearG}
Suppose hypothesis$(H')$ on $[0,T]$. Then there exists a linear map $\Gamma$ from the space of all $(\mathbb{P},\mathbb{F})$ local martingales into the space of c\`adl\`ag $\mathbb{G}$-predictable processes on $[0,T]$, with finite variation and null at the origin, such that, for any $(\mathbb{P},\mathbb{F})$ local martingale $X$, $\widetilde{X}:=X-\Gamma(X)$ is a $(\mathbb{P},\mathbb{G})$ local martingale on $[0,T]$. Moreover, if $X$ is a $(\mathbb{P},\mathbb{F})$ local martingale and $H$ is an $\mathbb{F}$ predictable process and $X$-integrable (in $(\mathbb{P},\mathbb{F})$), then $H$ is $\Gamma(X)$-integrable and $\Gamma(H\centerdot X)=H\centerdot \Gamma(X)$ on $[0,T]$. The operator $\Gamma$ will be called the drift operator.
\el

\

\subsection{Drift multiplier assumption}\label{viabilitydef}

In this paper we work especially with the drift operators having the following property.

\bd\label{assump1} (\textbf{Drift multiplier assumption}) We say that the drift multiplier assumption holds for the expansion $\mathbb{F}\subset \mathbb{G}$ on $[0,T]$ under $\mathbb{P}$, if 
\
\ebe 
\item
Hypothesis$(H')$ is satisfied for the expansion $\mathbb{F}\subset \mathbb{G}$ on the time horizon $[0,T]$ with a drift operator $\Gamma$;

\item
there exist $N=(N_1,\ldots,N_\mathsf{n})$ an $\mathsf{n}$-dimensional $(\mathbb{P},\mathbb{F})$ local martingale, and ${\varphi}$ an $\mathsf{n}$ dimensional $\mathbb{G}$ predictable process such that, for any $(\mathbb{P},\mathbb{F})$ local martingale $X$, $[N,X]^{\mathbb{F}\cdot p}$ exists, ${\varphi}$ is $[N,X]^{\mathbb{F}\cdot p}$-integrable, and $$
\Gamma(X)=\transp{\varphi}\centerdot [N,X]^{\mathbb{F}\cdot p}
$$
on the time horizon $[0,T]$. $N$ will be called the martingale factor and $\varphi$ will be called the integrant factor of the drift operator.
\dbe
\ed

An immediate consequence of the drift multiplier assumption is the following lemmas.

\bl\label{A-G-p}	
For any $\mathbb{F}$ adapted càdlàg process $A$ with $(\mathbb{P},\mathbb{F})$ locally integrable variation, we have $$
A^{\mathbb{G}\cdot p}=A^{\mathbb{F}\cdot p}+\Gamma(A-A^{\mathbb{F}\cdot p})=A^{\mathbb{F}\cdot p}+\transp\overline{\varphi}\centerdot[N,A-A^{\mathbb{F}\cdot p}]^{\mathbb{F}\cdot p}
$$
on $[0,T]$. In particular, for $R$ a $\mathbb{F}$ stopping time either $(\mathbb{P},\mathbb{F})$ predictable or $(\mathbb{P},\mathbb{F})$ totally inaccessible, for $\xi\in\mathbb{L}^1(\mathbb{P},\mathcal{F}_{R})$, $$
(\xi\ind_{[R,\infty)})^{\mathbb{G}\cdot p}
=
(\xi\ind_{[R,\infty)})^{\mathbb{F}\cdot p}+\transp\overline{\varphi}\centerdot(\Delta_{R}N\xi\ind_{[R,\infty)})^{\mathbb{F}\cdot p}
$$
on $[0,T]$.
\el

This lemma is proved in \cite{song-1505}.

\

\subsection{Martingale representation property}

The full viability will be studied under the martingale representation property. On the stochastic basis $(\Omega,\mathcal{A},\mathbb{P},\mathbb{F})$, consider a $d$-dimensional stochastic process $W$. We say that $W$ has the martingale representation property in the filtration $\mathbb{F}$ under the probability $\mathbb{P}$, if $W$ is a $(\mathbb{P},\mathbb{F})$ local martingale, and if all $(\mathbb{P},\mathbb{F})$ local martingale is a stochastic integral with respect to $W$. We say that the martingale representation property holds in the filtration $\mathbb{F}$ under the probability $\mathbb{P}$, if there exists a local martingale $W$ which possesses the martingale representation property. In this case we call $W$ a representation process.

\bethe\label{WY}
Supppose that $W$ has the martingale representation property in $(\mathbb{P},\mathbb{F})$. Then, the process $W$ satisfies the finite $\mathbb{F}$ predictable constraint condition. More precisely, there exist a finite number $\mathsf{n}$ of $d$-dimensional $\mathbb{F}$ predictable processes $\alpha_h, 1\leq h\leq \mathsf{n}$, such that$$
\Delta W 
= 
\sum_{h=1}^{\mathsf{n}}\alpha_h\ind_{\{\Delta W=\alpha_h\}}.
$$
\ethe

A very useful consequence of the finite predictable constraint condition is the following.

\bethe\label{boundedW}
If the martingale representation property holds in $\mathbb{F}$ under $\mathbb{P}$, there exists always a locally bounded representation process, which has pathwisely orthogonal components outside of a predictable thin set.
\ethe

The above two theorems are proved in \cite{song-mrp-drift}.

\

\subsection{Results on the full viability}

Here is an equivalence condition defined for expansion pair $\mathbb{F}\subset \mathbb{G}$.

\bcondition\label{1+fin}
For any $\mathbb{F}$ predictable stopping time $R$, for any positive random variable $\xi\in\mathcal{F}_R$, we have $\{\mathbb{E}[\xi|\mathcal{G}_{R-}]>0, R\leq T, R<\infty\}=\{\mathbb{E}[\xi|\mathcal{F}_{R-}]>0, R\leq T, R<\infty\}$. 
\econdition

\brem\label{rq:cR}
Clearly, if the random variable $\xi$ is already in $\mathcal{F}_{R-}$ (or if $\mathcal{F}_{R-}=\mathcal{F}_{R}$), the above set equality holds. Hence, a sufficient condition for Condition \ref{1+fin} to be satisfied is that the filtration $\mathbb{F}$ is quasi-left-continuous (cf. \cite[Definition 3.39]{HWY}).
\erem

The following theorem is proved in \cite{song-1505}. Let $T$ be a $\mathbb{G}$ stopping time.

\bethe\label{fullviability}
Suppose that $(\mathbb{P},\mathbb{F})$ satisfies the martingale representation property. Then, $\mathbb{G}$ is fully viable on $[0,T]$, if and only if the drift multiplier assumption \ref{assump1} (with the factors $N$ and $\varphi$) and Condition \ref{1+fin} are satisfied such that
\begin{equation}\label{fn-sur-fn}
\dcb
\transp{\varphi}(\centerdot[N^c,\transp N^c]){\varphi}\ \mbox{ is a finite process on $[0,T]$ and }\\
\\
\sqrt{\sum_{0<s\leq t\wedge T}\left(\frac{\transp {\varphi}_{s}\Delta_{s} N}{1+\transp{\varphi}_{s}\Delta_{s} N}
\right)^2},\ t\in\mathbb{R}_+,\
\mbox{ is $(\mathbb{P},\mathbb{G})$ locally integrable.}
\dce
\end{equation}

\ethe

\

\section{Recursive construction of multi-default time model}

In the previous sections, we have introduced in Definition \ref{fullviabilitydef} the notion of full viability of a general expansion pair $\mathbb{F}\subset \mathbb{G}$ and a result in Theorem \ref{fullviability} to ensure its validity. In present, we will consider more specifical expansions of $\mathbb{F}$ and study accordingly the full viability issue.

Fix a stochastic basis $(\Omega,\mathcal{A},\mathbb{P},\mathbb{F})$. Consider $n\in\mathbb{N}^*$ random times $\tau_1,\ldots,\tau_n$ defined on $(\Omega,\mathcal{A})$. We define $\mathbb{G}^{:0}=\mathbb{F}$, and, for $1\leq k\leq n$, define $\mathbb{G}^{:k}$ to be the progressive enlargement of the filtration $\mathbb{G}^{:(k-1)}$ with $\tau_k$ (cf. \cite[Chapitre IV]{J}). The filtration $\mathbb{G}^{:n}$ is an $n$-default times model.

\

\subsection{Full viability transmission}

We consider the full viability for the expansion pair $\mathbb{F}\subset \mathbb{G}^{:n}$. The following theorem shows that the full viability can pass harmonically through the recursive construction of the multi-default time models.

\bethe\label{viabilityrecursive}
For any $1\leq k\leq n$, let $T_{k}$ be a $\mathbb{G}^{:k}$ stopping time, with the property $T_{k}\leq T_{k-1}$ ($T_{0}=\infty$). Then, if the full viability holds for the expansion pair $\mathbb{G}^{:(k-1)} \subset \mathbb{G}^{:k}$ on the horizon $[0,T_{k}]$, $1\leq k\leq n$, the full viability holds for the expansion pair $\mathbb{F} \subset \mathbb{G}^{:n}$ on the horizon $[0,T_{n}]$. 
\ethe

\textbf{Proof.} The theorem is true for $n=1$. Suppose that the theorem is proved for $n=n_{0}, 1\leq n_{0}<n$. Let $S$ be an asset process in $\mathbb{F}$ possessing a deflator. Then, $S$ be an asset process in $\mathbb{G}^{:n_{0}}$ possessing also a deflator on the horizon $[0,T_{n_{0}}]$. As, by assumption, the full viability holds for the expansion pair $\mathbb{G}^{:n_{0}} \subset \mathbb{G}^{:n_{0}+1}$ on the horizon $[0,T_{n_{0}+1}] \subset [0,T_{n_{0}}]$, $S$ is again an asset process in $\mathbb{G}^{:n_{0}+1}$ possessing a deflator on the horizon $[0,T_{n_{0}+1}]$. \ \ok 

\

\subsection{Applicabiliy of Theorem \ref{fullviability} in recursive construction}

So, what is really essential is to establish the full viability for every expansion pair $\mathbb{G}^{:(k-1)}\subset \mathbb{G}^{:k}$. We will do this with Theorem \ref{fullviability}. This necessitates beforehand to ensure that every expansion pair $\mathbb{G}^{:(k-1)}\subset \mathbb{G}^{:k}$ satisfies the two conditions in section \ref{viabilitydef} and the martingale representation property. For this, we introduce three other conditions, which are easier to check when constructing multi-default time models. We consider the case of infinite horizon. 

\bassumption\label{exassump1}
\ebe
\item[.]
For any $1\leq k\leq n$, Hypothesis$(H')$ holds for the filtration expansion pair $\mathbb{G}^{:(k-1)}\subset\mathbb{G}^{:k}$ on the time horizon $[0,\infty]$ with a drift operator $\Gamma^{(k-1)|k}$.
\item[.]
The drift multiplier assumption holds for the expansion pair $\mathbb{G}^{:(k-1)}\subset\mathbb{G}^{:k}$ with the martingale factor $M^{(k-1)|k}$ (a $\mathbb{G}^{:(k-1)}$ local martingale) and the integrant factor $\overline{\psi}^{(k-1)|k}$ (a $\mathbb{G}^{:k}$ predictable process). 
\dbe
\eassumption

\bassumption\label{exassump2}
For any $1\leq k\leq n$, $\tau_k$ does not intercept the $\mathbb{G}^{:(k-1)}$ stopping times. 
\eassumption

\bassumption\label{exassump3}
For any $1\leq k\leq n$, the $s\!\mathcal{H}$ measure covering condition on $(0,\infty)$ (cf. \cite{JS}) is satisfied for the expansion pair $\mathbb{G}^{:(k-1)}\subset \mathbb{G}^{:k}$.
\eassumption

By Assumption \ref{exassump1}, the drift multiplier assumption for the expansion pair $\mathbb{G}^{:(k-1)}\subset \mathbb{G}^{:k}$ is ensured. To be able to apply Theorem \ref{fullviability}, it remains to establish the martingale representation property in every $\mathbb{G}^{:(k-1)}$ and also Assumption \ref{1+fin}.

\bethe\label{induct}
Suppose the assumptions \ref{exassump1} and \ref{exassump2} and \ref{exassump3}. Then, Hypothesis $(H')$ holds for the expansion pair $\mathbb{F}\subset \mathbb{G}^{:k}$. If $\mathbb{F}$ is quasi-left-continuous and satisfies the martingale representation property, for all $1\leq k\leq n$, $\mathbb{G}^{:k}$ also is quasi-left-continuous and satisfies the martingale representation property. 
\ethe

\textbf{Proof.} Consider the case where $k=1$. Hypothesis $(H')$ is satisfied for the expansion pair $\mathbb{G}^{:0}\subset \mathbb{G}^{:1}$ by Assumption \ref{exassump1}. 

The martingale representation property in $\mathbb{G}^{:1}$ is the consequence of Assumptions \ref{exassump2} and \ref{exassump3}, according to \cite[Theorem 5.1 and 5.2]{JS}.  

Note that, because of Assumption \ref{exassump2}, $\tau_1$ is $\mathbb{G}^{:1}$ totally inaccessible. By the drift multiplier Assumption \ref{exassump1} and the quasi-left-continuity in $\mathbb{G}^{:0}$, $\Gamma^{0|1}$ generates always continuous processes. This means, by Lemma \ref{A-G-p}, that the $\mathbb{G}^{:0}$ totally inaccessible stopping times have continuous compensators in $\mathbb{G}^{:1}$ so that they remain to be $\mathbb{G}^{:1}$ totally inaccessible stopping times. 

By the quasi-left-continuity, a $\mathbb{G}^{:0}$ local martingale can jump only at a $\mathbb{G}^{:0}$ totally inaccessible time. By \cite[Theorem 5.1]{JS}, the representation process in $\mathbb{G}^{:1}$ is linked with that in  $\mathbb{G}^{:0}$. Hence, the representation process in $\mathbb{G}^{:1}$ jumps only at $\mathbb{G}^{:1}$ totally inaccessible times. This implies that the filtration $\mathbb{G}^{:1}$ also is quasi-left-continuous. 

The theorem is proved for $k=1$. By induction, the theorem can be proved for all $1\leq k\leq n$.\ \ok

\

\subsection{Discussion on the martingale factor}

\bethe
Suppose the assumptions \ref{exassump1} and \ref{exassump2} and \ref{exassump3}. Suppose that $\mathbb{F}$ is quasi-left-continuous and satisfies the martingale representation property. Then, there exist an $\mathbb{F}$ local martingale $N^{(k-1)|k}$ and a $\mathbb{G}^{:k}$ predictable process $\overline{\varphi}^{(k-1)|k}$ such that, for any $\mathbb{F}$ local martingale $X$, $[N^{(k-1)|k},X]^{\mathbb{G}^{:(k-1)}\cdot p}$ exists, $\overline{\varphi}^{(k-1)|k}$ is $[N^{(k-1)|k},X]^{\mathbb{G}^{:(k-1)}\cdot p}$ integrable, and 
\begin{equation}\label{(k-1)N}
\Gamma^{(k-1)|k}(\widetilde{X}^{:(k-1)})=\transp\overline{\varphi}^{(k-1)|k}\centerdot [N^{(k-1)|k},X]^{\mathbb{G}^{:(k-1)}\cdot p},
\end{equation}
where $\widetilde{X}^{:(k-1)}$ denotes the martingale part of $X$ in $\mathbb{G}^{:(k-1)}$.
\ethe

\textbf{Proof.} Note that the initial formula is$$
\Gamma^{(k-1)|k}(\widetilde{X}^{:(k-1)})=\transp\overline{\psi}^{(k-1)|k}\centerdot [M^{(k-1)|k},\widetilde{X}^{:(k-1)}]^{\mathbb{G}^{:(k-1)}\cdot p},
$$
for a $\mathbb{G}^{:(k-1)}$ local martingale $M^{(k-1)|k}$. To prove the theorem, we only need to establish $$
 [M^{(k-1)|k},\widetilde{X}^{:(k-1)}]^{\mathbb{G}^{:(k-1)}\cdot p}<\!\!< [N^{(k-1)|k},X]^{\mathbb{G}^{:(k-1)}\cdot p},
$$
for some $\mathbb{F}$ local martingale $N^{(k-1)|k}$. Let $W$ be a locally bounded representation process for the martingale representation property in $\mathbb{F}$ (cf. Theorem \ref{boundedW}). By repeated applications of \cite[Theorem 5.1]{JS}, the representation process in $\mathbb{G}^{:(k-1)}$ is linked with $W$. With this link, because $X$ is an $\mathbb{F}$ local martingale, because of the continuity of the drift of $X$ in $\mathbb{G}^{:(k-1)}$ so that $\widetilde{X}^{:(k-1)}$ has jumps only at $\mathbb{F}$ stopping times, we can find a $\mathbb{G}^{:(k-1)}$ predictable process $H$ such that$$
[M^{(k-1)|k},\widetilde{X}^{:(k-1)}]^{\mathbb{G}^{:(k-1)}\cdot p}
=
\transp H\centerdot [\widetilde{W}^{:(k-1)},\widetilde{X}^{:(k-1)}]^{\mathbb{G}^{:(k-1)}\cdot p}
=
\transp H\centerdot [W,X]^{\mathbb{G}^{:(k-1)}\cdot p},
$$
where the last equality results from the continuity of the drifts of $W$ and of $X$ in $\mathbb{G}^{:(k-1)}$.\ \ok

\

\subsection{Recursive construction of the factors in the drift operators}

We can say more about the drift operator of the expansion pair $\mathbb{F}\subset \mathbb{G}^{:k}$ with the formula (\ref{(k-1)N}).

\bethe\label{induct2}
Suppose the assumptions \ref{exassump1} and \ref{exassump2} and \ref{exassump3}. Suppose that $\mathbb{F}$ is quasi-left-continuous and satisfies the martingale representation property. Then, Hypothesis $(H')$ is satisfied for the expansion pair $\mathbb{F}\subset \mathbb{G}^{:k}$, $1\leq k\leq n$, and the corresponding drift operator $\Gamma^{:k}$ satisfies the drift multiplier assumption whose factors $N^{:k},\overline{\varphi}^{:k}$ can be computed recursively by
$$
\dcb
^\top N^{:0}=0, \ ^\top \overline{\varphi}^{:0}=0,\\
^\top N^{:k}=(^\top N^{:(k-1)},\ ^\top N^{(k-1)|k}) \ \mbox{(forming a vector of one more dimension)},\\
^\top \overline{\varphi}^{:k} =(^\top\overline{\varphi}^{:(k-1)},\  \transp(1+\transp\overline{\varphi}^{:(k-1)} \gamma^{(k-1)|k} )\vdots\overline{\varphi}^{(k-1)|k}),\\

\dce
$$
where $\gamma^{(k-1)|k}=(\gamma^{(k-1)|k}_{i})_{1\leq i\leq d^{(k-1)|k}}$ is family of vectors which are determined by the relation $(\Delta N^{:(k-1)}{\centerdot}[N^{(k-1)|k}_i,\transp W])^{\mathbb{F}\cdot p} =\gamma^{(k-1)|k}_{i}{\centerdot}[N^{(k-1)|k}_i,\transp W]^{\mathbb{F}\cdot p}$, and $(1+\transp\overline{\varphi}^{:(k-1)} \gamma^{(k-1)|k} )$ is the vertical vector of components $(1+\transp\overline{\varphi}^{:(k-1)} \gamma^{(k-1)|k}_{i} )$, $1\leq i\leq d^{(k-1)|k}$ (the dimension of $N^{(k-1)|k}$).
\ethe

\textbf{Proof.} 
Note that $\gamma^{(k-1)|k}$ exists by \cite[Theorem 5.25 and Remark]{HWY}.
Note that $\Gamma^{:1}=\Gamma^{0|1}$ so that the drift multiplier assumption is satisfied for $k=1$. 

Suppose the induction assumption that the drift multiplier assumption is satisfied from $\mathbb{F}$ to $\mathbb{G}^{:(k-1)}$ with drift operator$$
\Gamma^{:(k-1)}(X)
=
\transp\overline{\varphi}^{:(k-1)}{\centerdot}[N^{:(k-1)},X]^{\mathbb{F}\cdot p}.
$$
Let us prove it for $k$.
Let $X$ be a $\mathbb{F}$ local martingale. Write the representation $X=\transp G\centerdot W$. The bracket process $[N^{(k-1)|k},X]^{\mathbb{F}\cdot p}$ is continuous, because of the quasi-left-continuity. By Lemma \ref{A-G-p}, for any component $N^{(k-1)|k}_i$,
$$
\dcb
&&
[N^{(k-1)|k}_i,X]^{\mathbb{G}^{:(k-1)}\cdot p}\\
&=&
[N^{(k-1)|k}_i,X]^{\mathbb{F}\cdot p}+\transp\overline{\varphi}^{:(k-1)}\centerdot[N^{:(k-1)},[N^{(k-1)|k}_i,X]-[N^{(k-1)|k}_i,X]^{\mathbb{F}\cdot p}]^{\mathbb{F}\cdot p}\\

&=&
[N^{(k-1)|k}_i,X]^{\mathbb{F}\cdot p}+\transp\overline{\varphi}^{:(k-1)}\centerdot[N^{:(k-1)},[N^{(k-1)|k}_i,X]]^{\mathbb{F}\cdot p}\\

&&\mbox{ because of the continuity $[N^{(k-1)|k}_i,X]^{\mathbb{F}\cdot p}$}, \\

&=&
[N^{(k-1)|k}_i,X]^{\mathbb{F}\cdot p}+\transp\overline{\varphi}^{:(k-1)}\centerdot(\Delta N^{:(k-1)}{\centerdot}[N^{(k-1)|k}_i,X])^{\mathbb{F}\cdot p}\\

&=&
[N^{(k-1)|k}_i,X]^{\mathbb{F}\cdot p}+\transp\overline{\varphi}^{:(k-1)}\centerdot(\Delta N^{:(k-1)}{\centerdot}[N^{(k-1)|k}_i,\transp W])^{\mathbb{F}\cdot p} G\\

&=&
[N^{(k-1)|k}_i,X]^{\mathbb{F}\cdot p}+\transp\overline{\varphi}^{:(k-1)}\gamma^{(k-1)|k}_{i}\centerdot [N^{(k-1)|k}_i,\transp W]^{\mathbb{F}\cdot p} G\\

&=&
[N^{(k-1)|k}_i,X]^{\mathbb{F}\cdot p}+\transp\overline{\varphi}^{:(k-1)}\gamma^{(k-1)|k}_{i}\centerdot [N^{(k-1)|k}_i, X]^{\mathbb{F}\cdot p} \\

&=&
(1+\transp\overline{\varphi}^{:(k-1)}\gamma^{(k-1)|k}_{i})\centerdot [N^{(k-1)|k}_i, X]^{\mathbb{F}\cdot p}.
\dce
$$
Let $\widetilde{X}^{:(k-1)}$ denote the martingale part of $X$ in $\mathbb{G}^{:(k-1)}$. We can write$$
\dcb
&&\Gamma^{(k-1)|k}(\widetilde{X}^{:(k-1)})
=
\transp\overline{\varphi}^{(k-1)|k}{\centerdot}[N^{(k-1)|k},\widetilde{X}^{:(k-1)}]^{\mathbb{G}^{:(k-1)}\cdot p}\\

&=&
\transp\overline{\varphi}^{(k-1)|k}{\centerdot}[N^{(k-1)|k},X]^{\mathbb{G}^{:(k-1)}\cdot p}\ \mbox{ because $[N^{:(k-1)},X]^{\mathbb{F}\cdot p}$ is continuous,}\\
&=&
\transp(1+\transp\overline{\varphi}^{:(k-1)}  \gamma^{(k-1)|k} )\vdots\overline{\varphi}^{(k-1)|k}\centerdot[N^{(k-1)|k},X]^{\mathbb{F}\cdot p}.
\dce
$$
Notice that $\widetilde{X}^{:k}=X-\Gamma^{:k}(X)$ and $\widetilde{X}^{:(k-1)}=X-\Gamma^{:(k-1)}(X)$ and also
$$
\dcb
&&\widetilde{X}^{:k} 
=
\widetilde{X}^{:(k-1)} - \Gamma^{(k-1)|k}(\widetilde{X}^{:(k-1)}) \\

&=&
X-\Gamma^{:(k-1)}(X) - \Gamma^{(k-1)|k}(\widetilde{X}^{:(k-1)}) \\
&=&
X-\transp\overline{\varphi}^{:(k-1)}{\centerdot}[N^{:(k-1)},X]^{\mathbb{F}\cdot p} - \transp(1+\transp\overline{\varphi}^{:(k-1)} \gamma^{(k-1)|k} )\vdots\overline{\varphi}^{(k-1)|k}\centerdot[N^{(k-1)|k},X]^{\mathbb{F}\cdot p}.
\dce
$$ 
The theorem is now proved with the formulas$$
\dcb
^\top N^{:k}=(^\top N^{:(k-1)},\ ^\top N^{(k-1)|k}),\\
^\top \overline{\varphi}^{:k} =(^\top\overline{\varphi}^{:(k-1)},\  \transp(1+\transp\overline{\varphi}^{:(k-1)} \gamma^{(k-1)|k} )\vdots\overline{\varphi}^{(k-1)|k}).\ \ok
\dce
$$

\

\section{Multi-default time $\natural$-model with full viability}\label{aboutdieze}

We continue to work on the stochastic basis $(\Omega,\mathcal{A},\mathbb{P},\mathbb{F})$. We want to construct a multi-default time model $\mathbb{G}^{:n}$ which satisfies the full viability. We want to do so with Theorem \ref{fullviability} and Theorem \ref{viabilityrecursive}. It is therefore necessary to construct beforehand a multi-default time model $\mathbb{G}^{:n}$, to which the two theorems are applicable. 

In the last section, Theorem \ref{induct} shows that, under Assumptions \ref{exassump1} and \ref{exassump2} and \ref{exassump3}, if the filtration $\mathbb{F}$ satisfies
the quasi-left-continuuity and
the martingale representation property, 
the filtrations $\mathbb{G}^{:k}$ preserve the same two properties, which make Theorem \ref{fullviability} and Theorem \ref{viabilityrecursive} applicable. Hence, what we have to do is to construct a multi-default time model satisfying Assumptions \ref{exassump1} and \ref{exassump2} and \ref{exassump3}.

\

\subsection{The $\natural$-model}

Can we have random times satisfying the above requirements ? We may mention honest times. However, \cite{FJS} has shown that a honest time model typically does not satisfy the no-arbitrage property on a horizon beyond the default time $\tau$. In contrast, the $\natural$-model introduced in \cite{JS2, songmodel} is perfect to meet the requirements.

Consider a stochastic basis $(\Omega,\mathcal{A},\mathbb{F},\mathbb{P})$. Let $L$ be a positive $\mathbb{F}$ locally bounded local martingale and $\Lambda$ be a non decreasing $\mathbb{F}$ adapted continuous process with $\Lambda_{0}=0$. We suppose that the non negative process $Z:=Le^{-\Lambda}$ satisfies $Z>0, 1-Z>0$ on $(0,\infty)$. Let $\mathbf{Y}$ be a multi-dimensional $\mathbb{F}$ locally bounded local martingale. We adopt the notion of smooth Markovian $\natural$-pair $(\mathbf{F}, \mathbf{Y})$ ($\mathbf{F}$ being a smooth Markovian coefficient) in \cite[subsection 3.3 and Theorem 3.7]{songmodel}. According to \cite[Theorem 4.8]{songmodel}, we have

\bethe
Under the above conditions, there exists a random time $\tau$ defined on (an isomorphic extension of) the probability space $(\Omega,\mathcal{A},\mathbb{P})$ such that 
\ebe
\item
$Z$ is the Azema supermartingale of $\tau$ with respect to $\mathbb{F}$, 
\item
the random time $\tau$ does not intercept the $\mathbb{F}$ stopping times,  
\item
Hypothesis $(H')$ holds for the expansion pair $\mathbb{F}\subset\mathbb{G}$, 
\item
the $s\!\mathcal{H}$ measure covering condition holds on $(0,\infty)$.
\dbe
Moreover, the drift operator takes the form
\begin{equation}\label{dies}
\Gamma(X)
=
\left(\ind_{(0,\tau]}\frac{1}{Z_-}
-
\ind_{(\tau,\infty)}\frac{1}{1-Z_{-}}
\right)\centerdot[M,X]^{\mathbb{F}-p}
+
\ind_{(\tau,\infty)}{\mathtt{p}}(\tau)\centerdot[\mathbf{Y},X]^{\mathbb{F}-p},
\end{equation}
where $M$ is the $\mathbb{F}$ martingale part of $Z$ ($M=e^{-\Lambda}\centerdot L$), $\mathtt{p}$ is a bounded $\mathbb{F}$ predictable process with parameter (coming form the derivative of $\mathbf{F}$). The pair $(\mathbb{F},\tau)$ will be called a $\natural$-model based on $(Z,\mathbf{F},\mathbf{Y})$.
\ethe

Note that the $s\!\mathcal{H}$ measure covering condition on $(0,\infty)$ is proved in \cite{JS}.  Note also that $L$ is assumed $\mathbb{F}$ locally bounded. As a corollary, we have the next results.

\

\subsection{A sufficient condition for the full viability}

Applying the previous theorem step by step in the construction of multi-default time models, we obtain the following corollary.

\bcor\label{diezeconstr}
Let, for $1\leq k\leq n$, $Z^{(k-1)|k}=L^{:(k-1)}e^{-\Lambda^{(k-1)|k}}$ be the Azema supermartingale of $\tau_{k}$ with respect to $\mathbb{G}^{:(k-1)}$. Suppose that $Z^{(k-1)|k}>0, 1-Z^{(k-1)|k}>0$, $L^{(k-1)|k}$ is locally bounded and $\Lambda^{(k-1)|k}$ is continuous. Suppose that $(\mathbb{G}^{:(k-1)},\tau_{k})$ is a smooth Markovian $\natural$-model based on $(Z^{(k-1)|k},\mathbf{F}^{(k-1)|k},\mathbf{Y}^{(k-1)|k})$ ($\mathbf{Y}^{(k-1)|k}$ being $\mathbb{G}^{:(k-1)}$ locally bounded). Then, Assumptions \ref{exassump1}, \ref{exassump2} and \ref{exassump3} are satisfied for the expansion pair $\mathbb{G}^{:k-1}\subset\mathbb{G}^{:k}$.
\ecor

Corollary \ref{diezeconstr} provides the conditions to apply Theorem \ref{fullviability} and Theorem \ref{viabilityrecursive}. We can now conclude.

\bethe\label{constr}
Suppose the same condition as in Corollary \ref{diezeconstr}. Suppose that 
\ebe
\item[$\bullet$]
$\frac{1}{Z^{(k-1)|k}}\ind_{(0,\tau_{k}]}$ and $\frac{1}{1-Z^{(k-1)|k}}\ind_{(\tau_{k},\infty)}$ are $\mathbb{G}^{:k}$ locally bounded processes and 

\item[$\bullet$]
$\mathbf{Y}^{(k-1)|k}$ is continuous. 
\dbe
Then, the full viability holds for the expansion pair $\mathbb{F} \subset \mathbb{G}^{:n}$, whenever $\mathbb{F}$ is quasi-left continuous and satisfies the martingale representation property.
\ethe

\brem
Notice that a multi-default time model may also be constructed with initial times (cf. \cite{pham}). But in general, a density process will not be as easy as $Z^{(k-1)|k}$ to be tackled. 
\erem

\

\end{document}